\RequirePackage{ifpdf}
\ifpdf 
\documentclass[pdftex]{sigma}
\else
\documentclass{sigma}
\fi

\begin{document}

\renewcommand{\PaperNumber}{036}

\FirstPageHeading

\renewcommand{\thefootnote}{$\star$}

\ShortArticleName{Conformally Covariant Operator}

\ArticleName{A Quartic Conformally Covariant Dif\/ferential\\ Operator for
Arbitrary Pseudo-Riemannian Manifolds (Summary)\footnote{This paper is a
contribution to the Proceedings of the 2007 Midwest
Geometry Conference in honor of Thomas~P.\ Branson. The full collection is available at
\href{http://www.emis.de/journals/SIGMA/MGC2007.html}{http://www.emis.de/journals/SIGMA/MGC2007.html}}}

\Author{Stephen M. PANEITZ}

\AuthorNameForHeading{S.M. Paneitz}

\Address{Deceased}

\ArticleDates{Received March 27, 2008 from Michael Eastwood; Published online March 30, 2008}

\Abstract{This is the original manuscript dated March $9^{\mathrm{th}}$~1983,
typeset by the Editors for the Proceedings of the Midwest Geometry Conference
2007 held in memory of Thomas Branson. Stephen Paneitz passed away on
September~$1^{\mathrm{st}}$~1983 while attending a conference in Clausthal and
the manuscript was never published. For more than 20 years these few pages were
circulated informally. In November 2004, as a service to the mathematical
community, Tom Branson added a scan of the manuscript to his website. Here we
make it available more formally. It is surely one of the most cited unpublished
articles. The dif\/ferential operator def\/ined in this article plays a key r\^ole
in conformal dif\/ferential geometry in dimension $4$ and is now known as the
Paneitz operator.}


\Classification{53A30; 58J70}

Let $M^\prime$ be a pseudo-Riemannian manifold of dimension $n>1$. The fourth
order operator def\/ined below (equation (\ref{three})) bears some analogy with
the second-order shift of the Laplace--Beltrami operator ${\mathcal{L}}$ by a
multiple of the scalar curvature~$R$, namely
\begin{equation}\label{one}{\mathcal{L}}+\frac{n-2}{4(n-1)}R\end{equation}
def\/ined by {\O}rsted, and studied by others in special cases. The relation
satisf\/ied by operators of the form (\ref{one}) is the following: if $M$ is
a manifold of dimension $n>1$ and $g_1$ and $g_2$ are two metric tensors on $M$
that satisfy $g_1=p^2g_2$ for some positive function~$p$, then
\begin{equation}
p^{(n+2)/2}\Big({\mathcal{L}}_{g_1}+\frac{n-2}{4(n-1)}R_{g_1}\Big)\phi
=\Big({\mathcal{L}}_{g_2}+\frac{n-2}{4(n-1)}R_{g_2}\Big)(p^{(n-2)/2}\phi)
\end{equation}
for all smooth functions~$\phi$, where ${\mathcal{L}}_{g_i}$ and $R_{g_i}$ are
the Laplace--Beltrami operator and scalar curvature determined by the
pseudo-Riemannian metric~$g_i$, for $i=1,2$. It is not necessary that the
scalar curvatures be constants. We will refer to the multipliers $p^{(n-2)/2}$
and $p^{(n+2)/2}$ as the ``initial'' and ``f\/inal'' multipliers, respectively.

Note that when $n=2$ (and only in this case), the initial multiplier is~$1$,
the multiplicative part of the operator (\ref{one}) is~$0$, the equation
${\mathcal{L}}_g\phi=0$ is independent of the metric~$g$, and thus these
solutions are simply related, modulo topology, to the harmonic $1$-forms
$\omega$ satisfying $d\omega=\delta\omega=0$, merely by exterior
dif\/ferentiation.

The quartic operator def\/ined below has an entirely comparable list of special
features in the physical case of dimension~$n=4$. Its solutions in this case
are related to the conformally invariant class of harmonic dif\/ferential forms,
which being forms of half the dimension of the manifold, are the two-forms; if
the metrics $g$ have Lorentzian signature the equations of harmonicity are the
Maxwell equations.

It seems reasonable to conjecture the existence of analogous operators of
arbitrary higher even order; indeed, powers of the f\/lat wave operator in
Minkowski space were shown to have a~conformal covariance property in
Jakobsen's thesis and work of Jakobsen--Vergne.

To def\/ine the quartic operator, we must assume that the dimension of the
manifold is neither $1$ nor~$2$. Given a pseudo-Riemannian metric $g$ on~$M$,
with associated Laplace--Beltrami operator~${\mathcal{L}}_g$, scalar
curvature~$R$, Ricci tensor components~$R_{ij}$ (tensorial notation referring
to coordinates $x^1,x^2,\dots,x^n$), metric components~$g_{ij}$ and covariant
derivatives $T_i=\nabla_{\frac{\partial}{\partial x^i}}$, def\/ine
\begin{gather}
{\mathcal{Q}}(g)\phi=({\mathcal{L}}_g)^2\phi+
T_j\left[\left(-\frac{4}{n-2}R^{ij}+\frac{n^2-4n+8}{2(n-1)(n-2)}g^{ij}R\right)
\frac{\partial}{\partial x^i}\phi\right]\nonumber\\
\phantom{{\mathcal{Q}}(g)\phi=}{} +\left[\frac{n-4}{4(n-1)}({\mathcal{L}}_gR)
-\frac{n-4}{(n-2)^2}R^{ij}R_{ij}+
\frac{(n-4)(n^3-4n^2+16n-16)}{16(n-1)^2(n-2)^2}R^2\right]\phi\label{three}
\end{gather}
for scalar functions~$\phi$ (the usual summation convention being followed).

\begin{theorem} Let $M$ be an arbitrary smooth manifold of
dimension $n>2$, and assume given two pseudo-Riemannian metrics $g_1$ and $g_2$
on $M$ related by $g_1=p^2g_2$ for a positive function~$p$. Then for all scalar
functions~$\phi$,
\begin{equation}\label{four}
p^{(n+4)/2}{\mathcal{Q}}(g_1)\phi={\mathcal{Q}}(g_2)(p^{(n-4)/2}\phi).
\end{equation}
When $n=4$,
\[{\mathcal{Q}}(g)\phi={\mathcal{L}}_g{}^2\phi
-2T_j\left[\left(R^{ij}-\frac{1}{3}g^{ij}R\right)\frac{\partial}{\partial x^i}\phi
\right],\]
and the kernel of this operator is independent of the metric~$g$, and also
invariant under the ordinary pointwise action of the conformal group $(M,g_1)$
(or equivalently $(M,g_2)$) on functions.
\end{theorem}

\begin{center}- - - - - - - - - - - - - - -\end{center}

In the special case of a conformally f\/lat four-dimensional
manifold $M$ of Lorentzian signature, we can sketch the relation to Maxwell's
equations alluded to above. For such an~$M$, any smooth local conformal
transformation extends uniquely to $\widetilde{M}={\mathbb R}\times S^3$. But now any
global distribution two-form solution of Maxwell's equations on $\widetilde{M}$
is automatically invariant under the central transformation of $\widetilde{M}$
sending $(t,p)$ into $(t+\pi,\mbox{ antipodal point to }p)$, as shown in the
work of I.E.~Segal. The associated compact quotient manifold
of~$\widetilde{M}$, denoted~$\overline{M}$, is casually identif\/iable with the
conformal compactif\/ication of four-dimensional Minkowski space. Now set
$g_c=dt^2-ds^2$ (the Einstein metric on ${\mathbb R}\times S^3$), and set
$X_0=\partial/\partial t$ and $\Delta$ the Laplacian on~$S^3$. Then
\[{\mathcal{Q}}(g_c)=(X_0{}^2-\Delta)^2+4X_0{}^2,\]
and
\[p^4{\mathcal{Q}}(g_c)={\mathcal{Q}}(g_f)=({\mathcal{L}}_{g_f})^2,\]
where $g_f=(dx^0)^2-(dx^1)^2-(dx^2)^2-(dx^3)^2$ is the f\/lat metric, and we have
$p^2g_f=g_c$. (Regarding~$S^3$ as~${\mathrm{SU}}(2)$, one can take
$p(t,u)=\frac{1}{2}\cos t+\frac{1}{4}{\mathrm{tr}}(u)$ ($t\in {\mathbb R}$,
$u\in{\mathrm{SU}}(2)$).)

\begin{theorem} ${\mathcal{Q}}(g_c)\phi=0$ has an
infinite-dimensional space of smooth solutions on $\overline{M}$ (all solutions
on $\widetilde{M}$ are lifted up from $\overline{M}$), and this equation
satisfies Huygens' principle. This space of solutions is the unique conformally
invariant closed subspace of the space of all smooth fields on $\overline{M}$
whose elements are determined by Cauchy data of a finite number of
$X_0$-derivatives on a~space-like surface.
\end{theorem}

Any two-form solution of Maxwell's equations on $\overline{M}$ is the exterior
derivative of a $1$-form $A$ on $\overline{M}$ satisfying $d\ast_cA=0$ and
$A(X_0)=0$ (gauge conditions; $\ast_c$ is the $\ast$-operator determined
by~$g_c$). The space of $1$-forms $A$ on $\overline{M}$ satisfying
$d\ast_c dA=0$ (Maxwell's equations) and $d\ast_cA=0$ and $A(X_0)=0$ is not
conformally invariant, but its direct sum with the conformally invariant space
of exterior derivatives of functions $\phi$ satisfying
${\mathcal{Q}}(g_c)\phi=0$ (or this space augmented by~$dt$, namely the
corresponding space of closed $1$-forms that are locally exterior derivatives
of functions $\phi$ satisfying ${\mathcal{Q}}(g_c)\phi=0$), is conformally
invariant.

Moreover, all the conformally invariant closed subspaces of the real smooth
$1$-forms on $\overline{M}$ satisfying Maxwell's equations are generated by the
above-noted conformally invariant subspaces and the obvious one of all closed
$1$-forms on~$\overline{M}$.

\LastPageEnding
\end{document}